\magnification=\magstep1
\input amstex
\define\R{\Bbb R}
\define\olo{\omega^{<\omega}}
\redefine\P{\text{\rm Power}}
\redefine\min{\text{\rm min}}
\redefine\S{\Bbb S}
\define\lr{L(\R)}
\define\oo{\omega^\omega}
\define\M{\Bbb M}
\define\sa{\S_\alpha}
\documentstyle{amsppt}
\pageheight{7.5in}
\NoBlackBoxes

\topmatter
\title
Isolating cardinal invariants
\endtitle
\author
Jind\v rich Zapletal
\endauthor
\affil
University of Florida
\endaffil
\abstract
There is an optimal way
of increasing certain cardinal invariants of the continuum.
\endabstract
\thanks
The author is partially supported by grants GA \v CR 201-00-1466 and NSF DMS-0071437.
\endthanks
\address
Department of Mathematics, University of Florida, Gainesville FL 32611
\endaddress
\email
zapletal\@math.ufl.edu
\endemail
\subjclass
03E17, 03E55, 03E60
\endsubjclass
\endtopmatter
\document
\head {0. Introduction}\endhead

The theory of cardinal invariants of the continuum is a large subfield of
set theory \cite {B2}. Its subject of study is the comparison of various cardinal
numbers typically defined as ``the smallest size of a set of reals with
certain properties''. Occasionally it is possible to prove inequalities
between these cardinals, but more often than not the inequalities are
independent of the usual axioms of set theory. Historically, certain
forcing extensions were identified as the standard tools for proving
these independence results; let me name various iterations of Sacks,
Cohen, Solovay or Laver real forcings as good examples. In this paper I prove that in a
certain precise sense some of these extensions are really the optimal tools for
establishing a broad syntactically defined class of independence results.
I will deal with the following class of invariants. 

\definition {0.1. Definition}
A {\it tame} invariant is one defined as $\min\{|A|:A\subset\R, 
\phi(A)\land\psi(A)\}$ where the quantifiers of $\phi(A)$ are restricted to the set $A$ or to the natural numbers and 
$\psi(A)$ is a sentence of the form
$\forall x\in\R\ \exists y\in A\ \theta(x, y)$ where $\theta$ is a formula whose quantifiers range over natural and 
real numbers only, without mentioning the set $A$. A real parameter is allowed in both formulas $\phi$ and $\psi$.
\enddefinition

Most cardinal invariants considered today are tame. For example:

\roster
\item"{$\circledast$}" $\frak a=\min\{|A|:A\subset [\omega]^\omega, \phi(A)\land\psi(A)\}$ where $\phi(A)=$``$A$ is an infinite set 
consisting of mutually almost disjoint sets'' and $\psi(A)=$``$\forall x\in [\omega]^\omega\ \exists y\in A\ x\cap y
$ is infinite''.
\item"{$\circledast$}" {\tt add}(meager)$=\min\{|A|:A\subset\R, \psi(A)\}$ where $\psi(A)=$``$\forall x\in\R\ \exists y\in A$ 
if $x$ codes a countable sequence of closed nowhere dense sets then $y$ codes a closed 
nowhere dense set not covered by their union''. In this case the sentence $\phi$ is not needed, that is we set $\phi=$
true.
\endroster

From these examples it is clear that in a definition of a tame invariant the sentence $\phi$ describes the 
internal structure of the set $A$ while $\psi$ is a statement about ``large size'' of the set $A.$ It is a routine matter
to write invariants like $\frak t, \frak u, \frak s$ as well as all the invariants in the Cichon
diagram \cite {B2} in a tame form. On the other hand, $\frak g$ and $\frak h$ apparently cannot be so written.

\proclaim {0.2. Theorem}
Suppose that there is a proper class of measurable Woodin cardinals. If $\frak x$ is a tame cardinal invariant
such that $\frak x<\frak c$ holds in some set forcing extension then $\frak x<\frak c$ holds in the iterated Sacks extension.
\endproclaim

Here, the iterated Sacks extension is obtained as usual by a countable support iteration of length $\frak c^+$
of Sacks forcing \cite {B2}. The theorem says that this extension is the optimal tool for proving the consistency
of inequalities of the type $\frak x<\frak c$ where $\frak x$ is a tame cardinal invariant. There are two
immediate consequences; I will state them in a rather imprecise form to retain their flavor. First, as in $P_{max}$
\cite {W2}, we get mutual consistency: if $\frak x_i:i\in I$ are tame invariants such that $\frak x_i
<\frak c$ is consistent for each $i\in I$ then even the conjunction of these inequalities is consistent.
 Restated, $\frak c$ cannot be
written as a nontrivial maximum of several tame invariants. And second, if $\frak x$ is a tame invariant
such that $\frak x<\frak c$ is consistent then so is $\aleph_1=\frak x<\frak c=\aleph_2$.

The proof of the theorem is flexible enough to give a host of related results.

\definition {0.3. Definition}
A cardinal invariant ${\frak y}$ can be {\it isolated} if there is a forcing $P_{\frak y}$ such that for every 
tame invariant $\frak x$,  if $\frak x<{\frak y}$ holds in some set forcing extension then it holds in the 
$P_{\frak y}$ extension.
\enddefinition

Thus the forcing $P_{\frak y}$ can be understood as increasing the invariant $\frak {y}$ in the gentlest way,
leaving all tame invariants smaller than $\frak y$ if possible. Hence the terminology. Theorem 0.2 says that
$\frak c$ can be isolated. I also have:

\proclaim {0.4. Theorem}
Suppose that there is a proper class of measurable Woodin cardinals. The following invariants can be isolated:

\roster
\item"{$\circledast$}" $\frak c$; $P_\frak c$ is the iterated Sacks forcing
\item"{$\circledast$}" $\frak b$; $P_\frak b$ is the iterated Laver forcing
\item"{$\circledast$}" $\frak d$; $P_\frak d$ is the iterated Miller forcing
\item"{$\circledast$}" $\frak h$; $P_\frak h$ is the iterated Mathias forcing
\item"{$\circledast$}" {\tt cov}(meager); one can use either a finite support or a countable support iteration of Cohen reals
\item"{$\circledast$}" {\tt cov}(null); use either a large measure algebra or a countable support iteration of Solovay reals
\item"{$\circledast$}" {\tt non}(strong measure zero); iterate forcings known as $PT_g$ \cite {B1}
\item"{$\circledast$}" {\tt add}(null); the forcing does not appear in published literature.
\endroster

\endproclaim

Amusingly enough, the proofs show that the minimum of any combination of invariants considered above 
can be isolated too by a countable support iteration in which the relevant forcings alternate. There are
invariants which cannot be isolated.  A good example is {\tt cof}(meager ideal)
since it can be written as max($\frak d$, {\tt non}(meager)). Both of the
inequalities $\frak d<${\tt cof}(meager), {\tt non}(meager)$<${\tt cof}(meager) are
consistent \cite {B2 2.2.11, 7.6.12, 7.5.8}. An invariant that cannot be isolated for a more complicated reason is
{\tt non}(meager). As shown in \cite {B2} {\tt cof}(meager)={\tt cov}($I_{ed}$) where $I_{ed}$
is the $\sigma$-ideal on $\oo$ generated by the sets $A_X=\{f\in\oo:\exists g\in X\ g\cap f$ is
infinite$\}$ as $X$ ranges over all countable subsets of $\oo.$ Now clearly {\tt cov}($I_{ed})\leq\sup\{
\frak d,$ {\tt cov}($I_{ed}(h)):h\in\oo\}$ where $I_{ed}(h)$ is the variation of the ideal $I_{ed}$ for the space
of all functions pointwise dominated by $h.$ However, the inequalities $\frak d<${\tt cov}$(I_{ed})$ as well
as {\tt cov}$(I_{ed}(h))<${\tt cov}$(I_{ed})$ for every fixed function $h\in\oo$ are consistent \cite {B2, S2}.
Ergo, the invariant {\tt non}(meager) cannot be isolated. This example was pointed out by Bartoszynski.

A curious twist of events occurs in the case of the tower number $\frak t$. 

\proclaim {0.5. Theorem}
Suppose that there is a proper class of measurable Woodin cardinals. There is a forcing $P_\frak t$ such
that for any tame invariant $\frak x$, if $\aleph_1=\frak x<\frak t$ holds in some forcing extension then it holds
in the $P_\frak t$ extension. 
\endproclaim

Thus it may be impossible to isolate $\frak t$ from invariants like $\frak p$ for which $\frak p<\frak t$ 
necessitates $\aleph_1<\frak p$. At the same time it is possible to choose the poset 
$P_{\frak t}$ to make $\frak t$ arbitrarily large.
Nothing like that occurs in the cases considered before. Also the forcing
$P_\frak t$ is undefinable, even though it is in some sense the expected thing.

The results stated above raise a number of obvious questions. For many invariants one would like to find 
out whether they can be isolated or not. If yes then what is the suitable forcing? If no, is there 
a clear reason? Above, I stated essentially everything I know in this direction at this point. That leaves
two of the invariants in the Cichon diagram without a status. Another issue is the use of large cardinal
hypotheses in the above theorems. Even though the proofs contain references to determinacy of certain integer
games of transfinite length and to $\Sigma^2_1$ absoluteness, I have no indication that the hypotheses used
are optimal or necessary at all.

The paper is organized as follows. The first section contains the analysis of the iterations 
of Sacks forcing from the descriptive set theoretic point of view. The complete proof of Theorem 0.2 can be 
found in the second section. In the third section I indicate the changes necessary to prove that 
$\frak b,\frak d, \frak h$ and {\tt non}(strong measure zero) can be isolated. The last section contains
the argument for Theorem 0.5.

The paper uses two important results whose proofs remain unpublished.

\proclaim {0.6. Fact}
($\Sigma^2_1$ absoluteness) (Woodin) Suppose that there is a proper class of measurable Woodin cardinals.
For every boldface $\Sigma^2_1$ sentence $\phi$, if $\phi$ holds in some generic extension then it holds
in every generic extension satisfying the continuum hypothesis.
\endproclaim

\proclaim {0.7. Fact}
(Transfinite projective determinacy) Suppose that there is a proper class of Woodin cardinals. Then every
integer valued game of every fixed transfinite countable length with projective outcome is determined. Moreover there is a 
winning strategy which is weakly homogeneous in every Woodin cardinal.
\endproclaim

The assumptions of the previous Fact are not optimal. Its proof consists of three parts. The determinacy of the
games was independently established by Neeman and Woodin. By a result of Martin \cite {Ma1} the games must have
winning strategies in a certain definability class. All sets in that definability class turn out to be weakly homogeneous
as shown by Neeman and Woodin independently.

The following fairly well known fact is the only property of weakly homogeneous sets we shall need.

\proclaim {0.8. Fact}
(Weakly homogeneous determinacy and absoluteness) \cite {W1} Suppose that $\delta$ is a supremum of Woodin 
cardinals with a measurable cardinal above it and
$T\subset(\omega\times Ord)^{<\omega}$ is a $<\delta$-weakly homogeneous tree. Then $\lr[p[T]]\models$AD and the theory
of the model $\lr[p[T]]$ with an arbitrary real parameter is invariant under forcing extensions of size $<\delta$.
\endproclaim

My notation follows the set theoretic standard set forth in \cite {J}, with one exception:
the concatenation of sequences $\vec r$ and $\vec s$ is denoted simply by $\vec r\vec s$.
Sequences of reals are denoted by $\vec r,\vec s\dots$ For a Polish space $X$ the expression Borel$(X)$ 
stands for the collection of all Borel subsets of $X.$ The spaces $\R^\alpha$ for a countable ordinal $\alpha$
are understood to come equipped with the product topology. A projective formula is one whose quantifiers range
over reals and integers only, and $\Sigma^2_1$ sentences are those of the form $\exists A\subset\R\ \theta(A)$
where $\theta$ is projective. Projective sets are usually confused with their definitions. For a tree $T$ the symbol $[T]$ stands for the set
of all its branches and $p[T]$ for the projection of this set into a suitable Polish space. AD$\R$ is the
statement ``all real games of length $\omega$ are determined''.  For a Woodin cardinal $\delta$
the expressions $\Bbb P_{<\delta}$ and $\Bbb Q_{<\delta}$ stand for the full nonstationary tower forcing on $\delta$
and its countably based variation respectively.
The reader is referred to \cite {J, B2, S1, W1} for all unfamiliar
concepts.

\head {1. The Sacks forcing}\endhead

The key to the proof of Theorem 0.2 is the understanding of Sacks forcing and its countable length countable support 
iterations in the context of determinacy. The well-known perfect set theorem can be restated to say that under ZF+AD
the Sacks forcing is (isomorphic to) a dense subset of the algebra $\P(\R)$ modulo the ideal of countable sets, ordered by inclusion.
It turns out that under the stronger determinacy hypothesis of ZF+DC+AD$\R$, for every countable ordinal $\alpha$
the countable support iteration of Sacks forcing of length $\alpha$ is a dense subset of the algebra $\P(\R^\alpha)$
modulo a suitable $\sigma$-ideal $I_\alpha$ on $\R^\alpha$. This is the driving idea behind the arguments.

\subhead {1.1. The geometric reformulation of Sacks forcing iterations}\endsubhead

First, it is necessary to restate the definition of the countable support iteration of countable length of Sacks forcing
in order to make the complexity analysis possible. A similar if not identical work was done by Kanovei in \cite {Ka}.
\definition {1.1.1. Definition}
For an ordinal $\alpha\in\omega_1$ define the poset $\Bbb S_\alpha$ to
consist of the nonempty Borel sets $p\subset\Bbb
R^\alpha$ satisfying these three conditions:
\roster
\item"{$\circledast$}" For every ordinal $\beta\in\alpha$ the set $p\restriction\beta=\{\vec s\in\Bbb
R^\beta:\exists \vec r\in p\ \vec s\subset\vec r\}$ is Borel. (The projection
condition; really for convenience only)
\item"{$\circledast$}" For every ordinal $\beta\in\alpha$ and every sequence
$\vec s\in p\restriction\beta$, the set
$\{t\in\Bbb R:\vec s\langle t\rangle\in p\restriction\beta+1\}$ is perfect. (The Sacks
condition)
\item"{$\circledast$}" For every increasing sequence $\beta_0\in\beta_1\in\dots$ of ordinals below $\alpha$
and every inclusion increasing sequence of sequences $\vec s_0\in p\restriction\beta_0, \vec
s_1\in p\restriction\beta_1\dots$, the sequence $\bigcup_n\vec s_n$ is in the set
$p\restriction\bigcup_n\beta_n$. (The countable support condition.)
\endroster
The sets $\Bbb S_\alpha$ are ordered by inclusion.
\enddefinition

It is not hard
to see that the posets
$\Bbb S_\alpha$ are naturally isomorphic to the countable support
iteration of
$\alpha$ many Sacks reals, if $\alpha\in\omega_1.$ If $G\subset\Bbb
S_\alpha$ is a generic filter then $G$ is given by the sequence $\vec
r_{gen}\in\Bbb R^\alpha$ with $\{\vec r_{gen}\}=\bigcap\{p:p\in G\}.$
This is done in the following lemma. The proof is completely unenlightening and should be skipped on the first reading of the paper.
It is however useful to notice that the argument depends only on the definability and properness of Sacks forcing.

\proclaim {1.1.2. Lemma}
Suppose $\alpha$ is a countable ordinal.
\roster
\item For every $\beta\in\alpha$ the poset $\Bbb S_\beta$ is naturally completely embedded
in $\sa$ and the factor $\sa/\Bbb S_\beta$ is naturally isomorphic to $\Bbb
S_{\alpha-\beta}$.
\item $\sa\Vdash$ for some unique sequence $\vec r_{gen}\in\Bbb R^\alpha$
the generic filter is just the set $\{p\in\check\Bbb S_\alpha:\vec r_{gen}\in p\}$.
\item If $\alpha$ is a limit ordinal then $\sa$ is the inverse limit of the posets $\{\Bbb
S_\beta:\beta\in\alpha\}.$
\item $\sa$ is proper.
\item For every countable elementary submodel $M$ of sufficiently large structure containing all the relevant information,
 for every forcing $P\in M$ adding a real $\dot s\in M$ and every
$P$-name $\dot p\in M$ for a condition in $\sa$ there is a Borel relation $B\subset\Bbb
R\times\Bbb R^\alpha$ so that
\itemitem{(a)} whenever $\beta\in\alpha$ then the relation
$B\restriction\beta=\{(s,\vec r)\in\Bbb R\times\Bbb R^\beta:\exists\vec t\
(s,\vec r\vec t)\in B\}$ is Borel
\itemitem{(b)} whenever
$(s,\vec r)\in B$ then
$s$ is
$M$-generic for
$P$,
$\vec r$ is
$M[s]$-generic for $\sa$ and $\vec r\in\dot p/s.$
\itemitem{(c)} for every $M$-generic real $s$ for $P$ the set $\{\vec r\in\Bbb
R^\alpha:(s,\vec r)\in B\}$ is a condition in $\sa$. It follows from (b) that this
condition strengthens $\dot p/s.$
\endroster
\endproclaim

Here (5) really amounts to saying that there is a constructive method for obtaining master
conditions. Note that by (2) and an absoluteness argument the condition obtained in (5c)
must be master for the model $M[s].$

\demo {Proof}
This is a completely standard simultaneous transfinite induction argument. I will show how (5) is obtained at
limit stages and why (1) holds at successor stages. The reader should refer to \cite
{S1,B2} for many similar arguments.
The following simple computation will be used throughout.

\proclaim {1.1.3. Claim}
Suppose $M$ is a countable transitive model of ZFC, $P\in M$ is a partial
order adding a single real $\dot s_{gen}$ and $\tau\in M$ is a $P$-name for a
real. Then
\roster
\item the set $A=\{s\in\Bbb R:$ the equation $s=\dot s_{gen}$ defines an
$M$-generic filter$\}$ is Borel
\item the function $\tau/s:A\to\Bbb R$ is a Borel function.
\endroster
\endproclaim

\demo {Proof}
I will prove (1); (2) is similar. The set $B$ of all $M$-generic filters on $RO(P)^M$ is
Borel in the product topology on Power$(RO(P)^M)$ since its elements $x$ are subject to the
Borel conditions ``$x$ is closed upwards", ``$x$ is a filter" and ``$x$ meets every open dense
set in $M$". The function $F:$ Power$(RO(P)^M)\to\Bbb R$ defined by $F(x)(n)=0$ if and only
if $\|\dot s_{gen}(\check n)=0\|^M\in x$ is continuous, and one to one on the set $B.$ Thus
the set $A$ is a one to one continuous image of a Borel set, therefore Borel. \qed
\enddemo

To see how Lemma 1.1.2(1) is obtained at a successor stage $\alpha=\beta+1$ we will prove that $\Bbb
S_\alpha\Vdash\vec r_{gen}(\beta)$ is $V[\vec r_{gen}\restriction\beta]$-generic Sacks
real. In order to do that, suppose $p_0\in\Bbb S_\alpha$ is an arbitrary condition and
$(q_0,\tau)\in\Bbb S_\beta*$Sacks is a condition such that $q\subset p\restriction\beta$ and
$q_0\Vdash\tau$ is a perfect subset of the set $\{t\in\Bbb R:\vec
r_{gen}\langle t\rangle\in p\}.$ It will be enough to produce a condition $p_1\subset p_0$ in
$\Bbb S_\alpha$ such that $p_1\Vdash \vec r_{gen}\restriction\beta\in q_0$ and $\vec
r_{gen}(\beta)\in\tau/\vec r_{gen}\restriction\beta.$ And indeed, if $M$ is a countable
elementary submodel of some large structure containing all relevant objects and
$q_1\subset q_0$ is a condition in $\Bbb S_\beta$ consisting of sequences $M$-generic for
this poset--and such a condition exists by the induction hypothesis (5)--, then we can put
$p_1=\{\vec r\in p_0:\vec r\restriction\beta\in q_1$ and $\vec
r(\beta)\in\tau/\vec r\restriction\beta\}$ and the condition $p_1\in\Bbb S_\alpha$ will be
as required.

Now suppose that $\alpha$ is a countable limit ordinal and Lemma 1.1.2(1)-(5) have been verified for
all ordinals $\beta\in\alpha.$ To prove (5) at $\alpha$ fix a countable elementary
submodel $M$ of some large structure $H_\lambda$ containing $\alpha$ and choose a forcing
$P\in M$ adding a real $\dot s$ and a $P$-name $\dot p\in M$ for a condition in $\sa.$
Choose an increasing sequence $\langle\alpha_n:n\in\omega\rangle$ of ordinals converging to
$\alpha$ starting with $\alpha_0=0$, and an enumeration $\langle\dot D_n:n\in\omega\rangle$
of all
$P$-names for open dense subsets of $\sa$ in $M.$ By induction on $n\in\omega$ choose
$P*\Bbb S_{\alpha_n}$-names
$\dot p_n\in M$ so that $\dot p=\dot p_0$ and 

\roster
\item"{$\circledast$}" $P*\Bbb S_{\alpha_n}\Vdash\dot p_n\in(\sa)^{V^P},\dot p_n\in\dot
D_{n-1},\vec r_{gen}\in \dot p_n\restriction\check\alpha_n$  
\item"{$\circledast$}" $P*\Bbb S_{\alpha_{n+1}}\Vdash$ if $\vec r_{gen}\in\dot
p_n\restriction\alpha_{n+1}$ then $\dot p_{n+1}\subset\dot p_n$.
\endroster

So for each $n\in\omega$ $\dot q_n=\{\vec t\in\Bbb R^{\alpha-\alpha_n}:\vec
r_{gen}\vec t\in\dot p_n\}$ is an $P*\Bbb S_{\alpha_n}$-name for a condition in
$\Bbb S_{\alpha-\alpha_n},$ $\dot q_n\in M$.

By the induction hypothesis there are Borel relations $B_n\subset\Bbb R\times\Bbb
R^{\alpha_n}\times\Bbb R^{\alpha_{n+1}-\alpha_n}$ so that

\roster
\item"{$\circledast$}" for every natural number $n$ and all $(s,\vec r,\vec t)\in B_n$ we have
that $s$ is an $M$-generic real for the poset $P,$ $\vec r$ is an $M[s]$-generic sequence
for $\Bbb S_{\alpha_n}$ and $\vec t$ is an $M[s][\vec r]$-generic sequence for $\Bbb
S_{\alpha_{n+1}-\alpha_n}$ such that $\vec t\in(\dot q_n/s,\vec r)\restriction
[\alpha_n,\alpha_{n+1})$
\item"{$\circledast$}" whenever $s$ is
an $M$-generic real for the poset $P$ and $\vec r$ is an $M[s]$-generic sequence for $\Bbb
S_{\alpha_n}$ then the set $\{\vec t\in\Bbb R^{\alpha_{n+1}-\alpha_n}:(s,\vec r,\vec t)\in
B_n\}$ is a condition in $\Bbb S_{\alpha_{n+1}-\alpha_n}.$
\endroster

Let $B\subset\Bbb R\times\Bbb R^\alpha$ be the relation given by $(s,\vec r)\in
B\leftrightarrow\forall n\in\omega\ (s,\vec r\restriction\alpha_n,\vec r\restriction
[\alpha_n,\alpha_{n+1})\in B_n.$ This is obviously a Borel relation, (5a) holds for it and
(5b, c) can be easily verified:

\roster
\item"{$\circledast$}" If $(s,\vec r)\in B$ then for all natural numbers $n<m$ we have $\vec
r\restriction\alpha_m\in(\dot p_n/s,\vec r\restriction\alpha_n)\restriction\alpha_m$ by the
choice of the names $\dot p_n,\dot q_n$ and the relations $B_n.$ By the countable support
condition applied to the sets $\dot p_n/s$ it must be the case that $\vec r\in\dot
p_n/s,\vec r\restriction\alpha_n$ for all $n\in\omega,$ in particular $\vec r\in\dot p/s$
and $\vec r$ is an $M[s]$-generic sequence for $\sa$.
\item"{$\circledast$}" Whenever $s$ is an $M$-generic real for the poset $P,$ we have $\{\vec
r\in\Bbb R^\alpha:(s,\vec r)\in B\}=\{\vec r\in\Bbb R^\alpha:\forall n\in\omega\ (s,\vec
r\restriction\alpha_n,\vec r\restriction[\alpha_n,\alpha_{n+1})\in B_n\}$ and the latter
set is easily verified to be a condition in $\sa.$
\endroster

Thus (5) has been proved for $\alpha.$ \qed
\enddemo

\subhead {1.2. The dichotomy}\endsubhead

The following is the key dichotomy and the only new result in this section.

\proclaim {1.2.1. Lemma}
\roster
\item (ZF+DC+AD$\R$) Suppose that $\alpha\in\omega_1$ and $A\subset\R^\alpha$. Then either there is a condition 
$p\in\Bbb S_\alpha$
with $p\subset A$ or there is a function $g:\R^{<\alpha}\to[\R]^{\aleph_0}$ such that $\forall\vec r\in A\ 
\exists\beta\in\alpha\ \vec r(\beta)\in g(\vec r\restriction\beta).$
\item (ZFC+there is a proper class of Woodin cardinals) Suppose that $\alpha\in\omega_1$ and $A\subset\R^\alpha$ is a projective 
set. Then the same dichotomy as in (1) holds for the set $A$.
\endroster
\endproclaim

The first item can be reworded thus: under ZF+DC+AD$\R$, the poset $\Bbb S_\alpha$ is a dense subset of the algebra $\P(\R^\alpha)$
modulo the $\sigma$-ideal $I_\alpha$ generated by the sets $B_g=\{\vec r\in\R^\alpha:\exists\beta\in\alpha\ \vec r(\beta)\in 
g(\vec r\restriction \beta)\}$ as $g$ varies through all functions from $\R^{<\alpha}$ to $[\R]^{\aleph_0}$. This is
a handsome way of putting things. However, the proof of (1) uses some hard unpublished theorems of Martin and Woodin
and works in a choiceless environment unfamiliar to some prospective readers. Since I will need the dichotomy for
projective sets only, I choose to include just the proof of (2). The assumption of (2) can be reduced to the existence of
$\omega_1$ Woodin cardinals.

\demo {Proof of Lemma 1.2.1(2)}
Let $\alpha\in\omega_1$ and $A\subset\R^\alpha$ be a projective set. Consider a real game of length $\alpha$
where players Adam and Eve play reals $s_\beta$ and $r_\beta$ respectively for $\beta\in\alpha$ so that
the real $s_\beta$ codes in some fixed way a countable set of reals and $r_\beta$ is not one of them.
Eve wins if the $\alpha$-sequence of her answers belongs to the set $A.$ Since real games of length
$\alpha$ are easily simulated by integer games of length $\omega\cdot\alpha$, by the Transfinite Determinacy Fact 0.7 the game
is determined and moreover there is a weakly homogeneous winning strategy. It is therefore enough to prove the following
two claims:

\proclaim {1.2.2. Claim}
Adam has a winning strategy iff there is a function $g:\R^{<\alpha}\to[\R]^{\aleph_0}$ such that $\forall\vec r\in A\ 
\exists\beta\in\alpha\ \vec r(\beta)\in g(\vec r\restriction\beta).$
\endproclaim

\proclaim {1.2.3. Claim}
Eve has a weakly homogeneous winning strategy iff there is a condition $p\in\S_\alpha$ with $p\subset A.$
\endproclaim

Now the first claim is a virtual triviality. The right-to-left direction of the second claim is not hard either. If $p\subset A$
for some condition $p\in\S_\alpha$ then Eve can defeat Adam merely making sure that at each stage $\beta\in\alpha$ the 
sequence $\vec r_\beta$ of answers she produced so far is in the set $p\restriction\beta$ and choosing her next answer from
the perfect set $\{t\in\R:\vec r_\beta\langle t\rangle\in p\restriction\beta+1\}$ minus the countable set coded
by Adam's challenge $s_\beta$. With a little care the choice can be made uniformly so that the winning strategy is not only
weakly homogeneous but Borel.

That leaves us with the left-to-right direction of Claim 1.2.3. Let $\sigma$ be a weakly homogeneous winning strategy for Eve.
Call a pair $\langle\vec s,\vec r\rangle$ of real sequences of length $\leq\alpha$ {\it good} if it represents a (partial)
play of the game in which Eve follows the strategy $\sigma.$ Thus there is a suitably weakly homogeneous tree $T$ whose
projection is the set of all good pairs of sequences of length $\leq\alpha$.

By transfinite induction on $\beta\leq\alpha$ prove that for every ordinal $\gamma\in\beta$ and every good pair $\langle \vec s_0,
\vec r_0\rangle\in\R^\gamma\times\R^\gamma$ there is a condition $p\in\S_{\beta-\gamma}$ such that for every sequence $\vec r_1\in p$
there is $\vec s_1\in\R^{\beta-\gamma}$ such that the pair $\langle \vec s_0\vec s_1,\vec r_0 \vec r_1\rangle\in
\R^\beta\times\R^\beta$ is good. This will clearly suffice considering the case $\beta=\alpha, \gamma=0$ and $\vec s_0=\vec r_0=0$ and
the fact that $\sigma$ is a winning strategy for Eve.

Suppose first that $\beta=\beta'+1$ is a successor ordinal and the induction hypothesis has been verified for $\beta'.$ Let 
$\gamma\in\beta$ and let $\langle \vec s_0,
\vec r_0\rangle\in\R^\gamma\times\R^\gamma$ be a good pair. By the induction hypothesis there is a condition $p'\in\S_
{\beta'-\gamma}$ such that for every $\vec r\in p'$ there is $\vec s$ such that the pair $\langle\vec s_0 \vec s,
\vec r_0\vec r\rangle$ is good. Now, whenever we have such a good pair then the set
$X_{\vec s,\vec r}=\{t\in\R:$ for some $u\in\R$ the pair $\langle\vec s_0 \vec s\langle u\rangle,
\vec r_0\vec r\langle t\rangle\rangle$ is good$\}$ must be uncountable; in the opposite case Adam would defeat
the strategy $\sigma$ by playing $\vec s_0,\vec s$ and then a code for the set $X_{\vec s,\vec r}.$ As $X_{\vec s,\vec r}\in
\lr[p[T]]$ and $\lr[p[T]]\models$AD the set $X_{\vec s,\vec t}$ must have a perfect subset. By the Weakly Homogeneous
Absoluteness 0.8, $p'\Vdash$ there is a sequence $\vec s$ such that $\langle\vec s_0\vec s,
\vec r_0\vec r_{gen}\rangle\in p[\check T]$ and there is a perfect set $c$ such that $\forall t\in c\exists u\in\R\ 
\langle\vec s_0\vec s\langle u\rangle,
\vec r_0\vec r_{gen}\langle t\rangle\rangle\in p[\check T]$. Pick $\S_{\beta'-\gamma}$ names $\vec s,\dot c$
for these two objects, let $M$ be a countable elementary submodel of a large enough structure containing all the relevant
information and using Lemma 1.1.2(5) find a condition $q\subset p'$ in $\S_{\beta'-\gamma}$ consisting of $M$-generic sequences only.
Then $p=\{\vec r\in\S_{\beta-\gamma}:\vec r\restriction\beta'\in q\land\vec r(\beta')\in\dot c/\vec r\restriction\beta\}$
is the sought condition in the poset $\S_{\beta-\gamma}.$ 

Now suppose $\beta\leq\alpha$ is a limit ordinal and the induction hypothesis has been verified up to $\beta.$ Suppose that
$\gamma\in\beta$ and $\langle\vec s_0,\vec r_0\rangle\in\R^\gamma\times\R^\gamma$ is a good pair. Let $\gamma=\beta_{-1}\in\beta_0\in
\beta_1\in\dots$ be an increasing $\omega$-sequence of ordinals converging to $\beta.$ By induction on $n\in\omega$ perform
the following three tasks:

\roster
\item"{$\circledast$}" Let $\vec r_{n+1}$ be the $\S_{\beta_n-\gamma}$-name for the part of the generic sequence of reals between $\beta_{n-1}$ and
$\beta_n$.
\item"{$\circledast$}" Choose an $\S_{\beta_n-\gamma}$-name $\vec s_{n+1}$ for a $\beta_n-\beta_{n-1}$-sequence of reals so that $\S_{\beta_n-\gamma}
\Vdash$ if there is a sequence $\vec s$ such that $\langle\vec s_0\vec s_1\dots \vec s_n\vec s, 
\vec r_0 \vec r_1 \dots\vec r_{n+1}\rangle\in p[\check T]$ then
$\vec s_{n+1}$ is such a sequence.
\item"{$\circledast$}" (Even for $n=-1$) Choose an $\S_{\beta_n-\gamma}$-name $\dot p_{n+1}$ for a condition in the forcing $\S_{\beta_{n+1}-\beta_n}$
such that $\S_{\beta_n-\gamma}
\Vdash$ if $\langle\vec s_0 \vec s_1 \dots \vec s_n \vec s_{n+1}, 
\vec r_0 \vec r_1 \dots\vec r_{n+1}\rangle\in p[\check T]$ then $\dot p_{n+1}\subset\{
\vec r\in\R^{\beta_{n+1}-\beta_n}:\exists \vec s\in\R^{\beta_{n+1}-\beta_n}\
\langle\vec s_0 \vec s_1 \dots 
\vec s_n \vec s_{n+1} \vec s, 
\vec r_0 \vec r_1 \dots\vec r_{n+1} \vec r\rangle\in p[\check T]\}$.
\endroster

For the third item note that the induction hypothesis has been proved up to $\beta$ and that by the Weakly Homogeneous Absoluteness
Fact 0.8 it holds up to $\beta$ even in the $\S_{\beta_n-\gamma}$ extension.

Now choose a countable elementary submodel $M$ of a large enough structure containing all the relevant information and use Lemma 1.1.2(5) to find Borel
relations $B_n\subset\R^{\beta_n-\gamma}\times\R^{\beta_{n+1}-\beta_n}$ for $n=-1,0,1,2\dots$ such that for all pairs $\langle\vec t_n,\vec t_{n+1}
\rangle\in B_n$ the sequence $\vec t_n$ is $M$-generic for $\S_{\beta_n-\gamma}$ and the sequence $\vec t_{n+1}\in\dot p_{n+1}/\vec t_n$ is
$M[\vec t_n]$-generic for $\S_{\beta_{n+1}-\beta_n}$, and moreover for every $M$-generic sequence $\vec t_n\in \R^{\beta_n-\gamma}$ the set
of all sequences $\vec t_{n+1}$ with $\langle\vec t_n,\vec t_{n+1}
\rangle\in B_n$ is a condition in the poset $\S_{\beta_{n+1}-\beta_n}.$ Let $p=\{\vec r\in\R^{\beta-\gamma}:$ for every $n=-1,0,1,\dots$ the
pair $\langle\vec r\restriction[\gamma,\beta_n),\vec r\restriction[\beta_n,\beta_{n+1})\rangle$ is in the relation $B_n\}.$ It is not difficult to
verify that $p\in\S_{\beta-\gamma}$ is the desired condition. \qed

\enddemo

\head {2. The absoluteness argument}\endhead

Towards the proof of Theorem 0.2, suppose that there is a proper class of measurable Woodin cardinals, $\frak x$
is a tame cardinal invariant, $\frak x=\min\{|A|:A\subset\R, \phi(A)\land\psi(A)\}$ where $\phi(A)$ is a statement 
quantifying over the natural numbers and elements of $A$, and $\psi(A)$ is a sentence of the form
$\forall x\in\R\ \exists y\in A\ \theta(x, y)$ where $\theta$ is a formula whose quantifiers range over natural and 
real numbers only, and suppose that $\frak x<\frak c$ holds in some set generic extension $V[G].$

Move into the model $V[G].$ There must be a set $A\subset\R$ such that $\phi(A)\land\psi(A)$ holds and $|A|<\frak c.$
I will prove that the Sacks forcing and its countable support iterations preserve the properties $\phi$ and $\psi$
of the set $A.$ Certainly $\phi(A)$ is preserved because of its simple syntactical form.
However the preservation of $\psi(A)$ could pose problems since some iteration 
$\S_\alpha$ could add a real $\dot x$ such that $\S_\alpha\Vdash\forall y\in\check A\ \lnot\theta(\dot x,y).$

\subhead {2.1. The countable case}\endsubhead

First consider the case of an arbitrary countable ordinal $\alpha\in\omega_1.$ Fix a condition $p\in\S_\alpha$
and an $\S_\alpha$-name $\dot x$ for a real. Strengthening the condition $p$ if necessary we may identify
$\dot x$ with a Borel function $\dot x:p\to\R$ with the understanding that the new real is the value
of this function on the generic $\alpha$-sequence of reals. I will show

$$\exists q\leq p\ \exists y\in A\ \forall \vec r\in q\ \theta(\dot x(\vec r),y).\tag *$$

Of course, then by projective absoluteness $q\Vdash\theta(\dot x,\check y)$ and as $p,\dot x$ were arbitrary,
$\S_\alpha\Vdash\forall x\in\R\ \exists y\in\check A\ \theta(x,y)=\psi(\check A)$ as desired. 

Suppose (*) fails. Then for every real $y\in A$ the set $B_y=\{\vec r\in p:\theta(\dot x(\vec r),y)\}$ contains no condition 
$q\leq p$ in the forcing $\S_\alpha$ as a subset. Since the sets $B_y$ are projective,
we can use the dichotomy 1.2.1 to find functions $g_y:\R^{<\alpha}\to [\R]^{\aleph_0}$ such that for every real $y\in A$
and every $\alpha$-sequence $\vec r\in p$ $\theta(\dot x(\vec r),y)$ implies $\exists\beta\in\alpha\ \vec r(\beta)\in g_y(\vec r\restriction
\beta).$ Now by transfinite induction on $\beta\in\alpha$ build an $\alpha$-sequence $\vec r\in p$ such that 
for every ordinal $\beta\in\alpha$ $\vec r\restriction\beta\in p\restriction\beta$ and $\forall y\in A\ \vec r(\beta)\notin g_y
(\vec r\restriction\beta)$. This is rather easy; at each level $\beta\in\alpha$ use the fact that $\bigcup_{y\in A}
g_y(\vec r\restriction\beta)$ is a set of size $|A|\cdot\aleph_0<\frak c$ while the set $\{t\in\R:(\vec r\restriction
\beta)\langle t\rangle\in p\restriction\beta+1\}$ is perfect, therefore of size $\frak c$ and so must contain
a real not in the above union. Now look at the real $\dot x(\vec r).$ By the choice of the functions $g_y$ and the sequence
$z$ we should have $\forall y\in A\ \lnot\theta(\dot x(\vec r),y),$ contradicting the property $\psi$ of $A$. (*) follows.

\subhead {2.2. The uncountable case}\endsubhead

The results of the previous subsection can be extended by a rather standard argument to show that
for every ordinal $\alpha$ the countable support iteration $\sa$ of Sacks forcing of length $\alpha$ preserves the statement $\psi(A)$.
Just use the following lemma:

\proclaim {Lemma 2.2.1}
(ZFC+projective absoluteness) Suppose that $\theta(x,y)$ is a projective formula and $A\subset\R$ 
is a set such that for every ordinal $\beta\in\omega_1,$
every condition $p\in\S_\beta$ and every Borel function $f:p\to\R$ there is a condition $q\leq p$ and a real $y\in A$
such that for every sequence $\vec r\in q,$ $\theta(f(\vec r),y)$ holds. Then for every ordinal $\alpha,$ $\sa\Vdash
\forall x\in\R\exists y\in\check A\ \theta(x, y).$
\endproclaim

Note that the assumptions of the lemma were shown to hold in the model $V[G]$ in the previous subsection.

\demo {Proof}
First, a small observation. Suppose $\beta\in\omega_1$ and $\alpha$ are
ordinals and $\pi:\beta\to\alpha$ is an increasing function. Then $\pi$
can be naturally extended into an order-preserving map $\pi:\Bbb
S_\beta\to\Bbb S_\alpha$ where $\pi(p)$ is the unique condition in $\Bbb
S_\alpha$ with support $\pi''\beta$ such that $\forall\gamma\in\beta\
\pi(p)\restriction\pi(\gamma)\Vdash_{\Bbb
S_{\pi(\gamma)}}(\pi(p))(\pi(\gamma))=\{t\in\Bbb R:\langle\dot
r_{\pi(\xi)}:\xi\in\gamma\rangle\langle t\rangle\in p\restriction\gamma+1\}$, where $\dot
r_\zeta$ is the
$\zeta$-th Sacks generic real. It is not hard to see that $\pi(p)\Vdash_{\Bbb
S_\alpha}\langle\dot r_{\pi(\xi)}:\xi\in\beta\rangle\in p.$

Now suppose that $\theta,A$ satisfy the assumptions of the lemma, $\alpha$ is an ordinal,
$q_0\in\S_\alpha$ is a condition and $\dot x$ is an $\sa$-name for a real. I will produce a condition $q_1\leq q_0$ and a real
$y\in A$ such that $q_1\Vdash\theta(\dot x,\check y).$ This will
prove the lemma. Choose a countable elementary submodel $M$ of some large
structure containing all relevant objects and let
$\beta=o.t.M\cap\alpha$ and $\pi:\beta\to\alpha$ be the inverse of the
transitive collapse. A standard countable support iteration argument similar to the proof
of Lemma 2.2(5) gives a condition $p_0\in\Bbb S_\beta$ such that $\pi(p_0)\leq q_0$ and
for every $\vec r\in p_0$ the sequence $\vec r\circ\pi^{-1}$ is
$M$-generic for the poset $\Bbb S_\alpha$. Let $f:p_0\to\Bbb R$ be the
Borel function defined by $f(\vec r)=\dot x/\vec r\circ\pi^{-1}.$ Thus
$\pi(p_0)\Vdash\dot x=\dot f(\langle\dot
r_{\pi(\xi)}:\xi\in\beta\rangle)$. The assumptions of the lemma can be now employed to provide a real $y\in A$
 and a condition $p_1\leq p_0$ such that
$\forall\vec r\in p_1\ \theta(f(\vec r),y)).$ By the projective absoluteness and the last
sentence of the first paragraph of this proof, setting
$q_1=\pi(p_1)$ we have $q_1\leq q_0,$ $q_1\Vdash_{\Bbb
S_\alpha}\theta(\dot x,\check y)$ as desired. \qed
\enddemo

\subhead {2.3. The wrap-up}\endsubhead

To restate the above work, let

$$\chi(A)=\forall\alpha\in\omega_1\ \forall p\in\S_\alpha\ \forall\dot x:p\mapsto\R\ {\text {\rm Borel}}\ \exists
y\in A\ \exists q\in\S_\alpha\ q\leq p\land\forall\vec r\in q\ \theta(\dot x(\vec r),y)$$

\noindent Note that $\chi(A)$ is a projective statement about the set $A\subset\R.$ We proved that $V[G]\models\chi(A)$ and that
$\chi(A)$ implies in ZFC+projective absoluteness that for every ordinal $\alpha$ $\S_\alpha\Vdash\psi(\check A)$.
Again,

$$V[G]\models\exists A\subset\R\ \phi(A)\land\chi(A).$$

\noindent Note that the sentence on the right hand side of the $\models$ sign is $\Sigma^2_1.$

Now back to the ground model $V.$ Suppose first that $V$ satisfies the continuum hypothesis. Then by the $\Sigma^2_1$ 
Absoluteness Fact 0.6, $V\models\exists A\subset\R\ \phi(A)\land\chi(A).$ Fix a set $A\subset\R$ with $\phi(A)\land
\chi(A)$ and iterate Sacks reals $\omega_2$ times with countable support to get a model $V[H].$ By the above work,

$$V[H]\models\phi(A)\land\psi(A), \frak x\leq|A|\leq|\frak c^V|=\aleph_1<\frak c={\aleph_2}^V=\aleph_2$$

\noindent as desired. If the continuum hypothesis fails in the ground model $V,$ iterate the Sacks reals $\frak c^+$ many times
anyway to get the model $V[H].$ Let $V[K]\subset V[H]$ be the intermediate extension given by the first
$\omega_1$ many generic reals. As is well known, $V[K]\models CH$ and $V[H]$ is an $\omega_2$ iterated Sacks extension
of the model $V[K].$ One can then repeat the above argument with $V$ replaced with $V[K]$ to see that $V[H]\models
\frak x<\frak c.$ Theorem 0.2 follows.

\head {3. Other invariants}\endhead

Many invariants of the form {\tt cov}($I$), where $I$ is a Borel generated $\sigma$-ideal on the real line,
can be isolated by a countable iteration of the forcing Borel($\R$) modulo $I,$ and the proof follows closely the scenario
of the previous two sections. It is just sufficient
to verify that this forcing is proper (this fact is used in setting up the geometric representation of the iteration as in
Subsection 1.1), that under AD every $I$-positive set of reals has a Borel positive subset (this is needed for the
successor step in the proof of Claim 1.2.3) and that {\tt cov}($I$)={\tt cov}($I\restriction B$) for every positive
Borel set $B,$ this is tacitly used in the proof of (*) in Subsection 2.1. The invariants $\frak c=${\tt cov}(countable),
$\frak b,\frak d$ and some others conform exactly to this scenario. For the invariants {\tt non}(strong measure zero)
and $\frak h$ further changes are necessary. 

\subhead {3.1. The dominating number}\endsubhead

Clearly the dominating number is the covering number of the ideal of bounded subsets of $\oo$. The countable support
iteration of Miller forcing \cite {Mi2} will isolate it as the following two lemmas show.

\proclaim {3.1.1. Lemma} \cite {Ke}
 Every Borel unbounded subset of $\oo$ contains all branches of some superperfect tree.
Under AD this generalizes to all unbounded sets.
\endproclaim

\noindent Thus the Miller forcing is a dense subset of the factor
algebra Borel
$(\oo)$ modulo the bounded sets. 

\proclaim {3.1.2. Lemma}
For every superperfect tree $T\subset\olo$ there is a continuous function $F:\oo\to[T]$ such
that preimages of bounded sets are bounded.
\endproclaim

\demo {Proof}
Thinning the tree $T$ out if necessary we may assume that every splitnode of $T$ has in fact
infinitely many immediate successors. The natural homeomorphism $F:\oo\to[T]$ will have the
required property. \qed
\enddemo

\noindent Thus $\frak d=${\tt cov}(bounded)={\tt cov}(bounded ideal restricted to $B$) for every Borel unbounded set $B\subset\oo.$
The argument in Sections 1 and 2 now goes through with the obvious changes, replacing $\frak c$ with $\frak d$, 
the countable ideal with the bounded ideal,
and the Sacks condition in Definition 1.1.1 with the obvious Miller condition.

\subhead{3.2. The bounding number}\endsubhead

The countable support iteration of Laver reals \cite {L} isolates $\frak
b$. Consider the
$\sigma$-ideal $I_L$ on $\oo$ generated by the sets $A_g=\{f\in\oo:$ for infinitely many
$n\in\omega$
$f(n)\in g(f\restriction n)\}$ where $g$ varies through all functions from $\olo$ to
$\omega$. We have the almost obvious

\proclaim {3.2.1. Lemma}
$\frak b=${\tt cov}$(I_L).$
\endproclaim

\demo {Proof}
The map $G:\oo\to I_L$ defined by $G(f)=A_g$ where $g(t)=f(|t|)$ for every sequence
$t\in\olo,$ has the property that preimages of non-covering subsets of $I_L$ are bounded.
This proves that {\tt cov}$(I_L)\leq \frak b$. On the other hand, fixing an enumeration
$\{u_n:n\in\omega\}$ of $\olo$,  the map
$H:I_L\to\oo$ sending the set $A_g$ to the function $f:n\mapsto g(u_n)$, has the property
that preimages of bounded sets do not cover the whole real line. Thus $\frak
b\leq${\tt cov}$(I_L).$ \qed
\enddemo

As in the previous subsection, I will prove that the Laver forcing is a dense subset of the algebra Borel$(\oo)$ modulo the ideal $I_L$:

\proclaim {3.2.2. Lemma}
 Every Borel $I_L$-positive set contains all branches of some Laver tree. Under AD this
generalizes to all $I_L$-positive sets.
\endproclaim

\demo {Proof}
Suppose $A\subset\oo$ is a set and define an infinite game by letting players Adam and Eve play
sequences $t_n\in\olo$ and bits $b_n\in 2$ respectively, observing the following rules:
$b_0=1$ and whenever Eve accepts a sequence $t_n$--that is, plays $b_n=1$--then
Adam submits one-step extensions $t_{n+1}, t_{n+2}, \dots$ of $t_n$ until Eve
accepts one of them. The last number on the sequences $t_{n+1}, t_{n+2},\dots$ must
increase. Adam wins if either Eve accepted only finitely many times or else
$\bigcup\{ t_n:b_n=1\}\in A.$ The following two claims will complete the proof of the lemma
\cite {Ma2}:

\proclaim {3.2.3. Claim}
Adam has a winning strategy if and only if the set $A$ contains all branches of some
Laver tree.
\endproclaim

\demo {Proof}
For the right to left direction fix a Laver tree $T$ with $[T]\subset A.$ Let Adam set
$t_0=$trunk of $T$, and if $t_n\in T$ has been played and accepted by Eve then let
Adam submit immediate successors of the node $t_n$ in the tree $T$ in the increasing
order until Eve accepts one of them. This is obviously a winning strategy for Adam.

For the left to right direction let $\sigma$ be a winning strategy for Adam and let
$T\subset\oo$ be the tree of all sequences that can possibly arise in a run of the game
$G_A$ in which Adam follows the strategy $\sigma$. Note that for each node $t\in T$
there is a unique shortest run $\tau(t)$ such that it respects $\sigma$ and $t$ occurs in
it, and if
$t\subset s$ are both in the tree $T$ then $\tau(t)\subset\tau(s).$ It follows that every
branch $f\in [T]$ is a result of the run $\bigcup\{\tau(t):t\subset f\}$ and therefore must
belong to the set
$A.$ It is also clear from the definition of the game $G_A$ that $T$ is a Laver tree with
trunk $\sigma(0)$. \qed
\enddemo

\proclaim {3.2.4. Claim}
Eve has a winning strategy if and only if $A\subset A_g$ for some function
$g:\olo\to\omega.$
\endproclaim

\demo {Proof}
For the right to left direction fix a function $g$ such that $A\subset A_g.$ Let Eve
accept a sequence $t_n,$ a one-step extension of some previously accepted sequence $t_m$ as
soon as the last number on $t_n$ exceeds $g(t_m).$ The result of such a play must fall
outside of the set $A_g$ and therefore this is a winning strategy for Eve.

For the left to right direction let $\sigma$ be a winning strategy for Eve. For every
sequence $s\in\olo$ let $T_s$ be the tree of all sequences that can be accepted by Eve
in some run of the game where he follows the strategy $\sigma$ and Adam plays
$t_0=s.$ It follows that for all sequences $s\subset t,$ if $t\in T_s$ then all but
finitely many one-step extensions of $t$ must belong to the tree $T_s$--otherwise Adam
could win by first getting to $t$ and then submitting all the one-step extensions of $t$
which do not belong to the tree $T_s.$ Also, $[T_s]\cap A=0$ for all $s\in\olo.$ To see
this, fix a branch $f\in [T_s]$ and define $S$ to be the tree of all partial runs of the
game $G_A$ in which Adam set $t_0=s,$ Eve followed the strategy $\sigma$ and the
last move of Adam was accepted and it is an initial segment of the branch $f.$ The tree
$S$ is ordered by extension. It follows from the ``increasing'' rule of the game $G_A$ that
the tree $S$ is finitely branching--each run $\tau\in S$ has at most $2^{f(n)}$ immediate
successors where $n$ is the length of the last move of $\tau$. Also, the tree $S$ has
height $\omega$, so it must be illfounded. Any infinite branch of the tree $S$ yields a run
of the game $G_A$ following the winning strategy $\sigma$ whose result was the function
$f.$ Thus $f\notin A.$

Now define a function $g:\olo\to\omega$ by setting $g(t)=$ an integer such that for every
$s\subset t,$ if $t\in T_s$ then $g(t)$ is larger than all of the finitely many numbers $n$
such that $t\langle n\rangle\notin T_s.$ I claim that $A\subset A_g.$ If this were not true
then there would be a function $f\in A$ such that for some $n\in\omega,$ for all larger
numbers $m$ necessarily $g(f\restriction m)\in f(m).$ But then $f\in [T_{f\restriction
n}]$ by the definition of the function $g,$ so by the previous paragraph $f\notin A.$ A
contradiction! \qed
\enddemo
\enddemo 

The last thing that must be verified before unleashing the technology developed in Sections 1 and 2
is that {\tt cov}$(I_L)=${\tt cov}$(I_L\restriction B)$ for every Borel $I_L$ positive set $B\subset\oo$:

\proclaim {3.2.5. Lemma}
For every Laver tree $T$ there is a continuous function $F:\oo\to[T]$ such that preimages
of $I_L$-small sets are $I_L$-small.
\endproclaim

\demo {Proof}
The natural homeomorphism $F:\oo\to[T]$ has the required property. \qed
\enddemo

\subhead {3.3. The uniformity of the strong measure ideal}\endsubhead

This invariant has a definition that is not suitable for our purposes for syntactical reasons.
I will use the following
combinatorial characterization of this invariant. For a function $g\in\oo$ let $I_{ie}(g)$
be the $\sigma$-ideal on $\Pi_ng(n)$ generated by the sets $A_f=\{h\in\Pi_ng(n):h\cap f$ is
finite$\}.$ Then

\proclaim {3.3.1. Lemma} \cite {B2 8.1.14, M1}
{\tt non}(strong measure zero)$=\min\{${\tt cov}$(I_{ie}(g)):g\in\oo\}.$
\endproclaim

Thus a natural attempt at isolating {\tt non}(strong measure zero) is the countable support iteration of the forcings
Borel($\Pi_ng(n)$) modulo the ideal $I_{ie}(g)$ for all possible (names for) functions $g\in\oo$. The following
two lemmas show that this attempt will actually work. Lemma 3.3.2 gives us the representation of the forcings suitable
to prove that they satisfy Axiom A, and yields the crucial dichotomy. Lemma 3.3.5 provides the necessary homogeneity
in the covering number.

Fix a function $g\in\oo$. A nonempty tree $T\subset\olo$ will be called $g${\it -thick} if
the sequences in $T$ are everywhere dominated by the function $g$, and for every
sequence $t\in T$ there is a natural number $n$ such that for every $m\in g(n)$ there is an
extension $s\in T$ of the sequence $t$ such that $s(n)=m.$ It is quite obvious that if $T$
is a $g$-thick tree then
$[T]\subset\Pi_n g(n)$ is an $I_{ie}(g)$-positive set. In fact,

\proclaim {3.3.2. Lemma}
For every function $g\in\oo,$ every Borel $I_{ie}(g)$-positive set contains all branches of
some $g$-thick tree. Under AD this generalizes to all $I_{ie}(g)$-positive sets.
\endproclaim

\demo {Proof}
Let $g\in\oo$ be a function and let $A\subset\Pi_ng(n)$ be a set. Define a game $G_A$ by
setting 

\smallskip
\settabs 8\columns
\+Adam&$t_0 ,n_0$&&$t_1, n_1$&&$t_2, n_2$&&\dots\cr
\+Eve&&$m_0$&&$m_1$&&$m_2$&\dots\cr
\smallskip
\noindent where $n_0, n_1, \dots$ is an increasing sequence of natural numbers, $m_i\in
g(n_i)$ and $0=t_0\subset t_1\subset\dots$ are sequences of natural numbers dominated by the
function $g$, $dom(t_i)\in n_i$ and $t_{i+1}(n_i)=m_i.$ Adam wins if $\bigcup t_n\in A.$
The following two claims will complete the proof of the lemma \cite {Ma2}:

\proclaim {3.3.3. Claim}
Adam has a winning strategy if and only if the set $A$ contains all branches of some
$g$-thick tree.
\endproclaim

\demo {Proof}
For the right to left direction fix a $g$-thick tree $T$ with $[T]\subset A.$ Adam will
easily win by making sure that for each of his moves $t_i\in T$, and that $n_i$ is such that
for every $m\in g(n_i)$ there is an extension $s\in T$ of the sequence $t_i$ such that
$s(n_i)=m_i.$

For the left to right direction fix a winning strategy $\sigma$ for Adam. Let $T$ be
the closure under initial segment of the set of all sequences arising in partial runs of
the game $G_A$ in which Adam follows the strategy $\sigma$. It is immediately clear
that $T$ is a $g$-thick tree and if $h$ is a branch through $T$ then there is a unique run
of the game in which Adam follows the strategy $\sigma$ and obtains the function $h.$
Ergo, $[T]\subset A$. \qed
\enddemo

\proclaim {3.3.4. Claim}
Eve has a winning strategy if and only if $A\subset\bigcup_k A_{f_k}$ for some
functions $f_k\in\Pi_n g(n),k\in\omega.$
\endproclaim

\demo {Proof}
For the right to left direction let $A\subset\bigcup_k A_{f_k}.$ Eve will easily win
by fixing a bookkeeping function $b:\omega\to\omega$ such that for every number $k$ the set
$b^{-1}\{k\}$ is infinite, and then playing $m_i=g_{b(i)}(n_i).$

For the left to right direction let $\sigma$ be a winning strategy for Eve. For each
partial run $\tau$ of the game $G_A$ where Eve followed the strategy and Adam
made the last move $t_i$ let $f_\tau\in\Pi_n g(n)$ be the
function defined by $f_\tau(n)=\sigma(\tau\langle n\rangle)$. Then necessarily
$A\subset\bigcup_\tau A_{f_\tau}.$ If this failed, then there would be a function $h\in A$
with infinite intersection with each $f_\tau.$ And then Adam could beat the strategy
$\sigma$ by inductively constructing a run of the game which respects the strategy
$\sigma$ and results in the function
$h.$ Assuming that the partial run $\tau_i$ has been constructed so that Adam made a
last move $t_i\subset h$ in it, he finds a number $n_i$ such that $h(n_i)=f_{\tau_i}(n_i)$
and the game continues into $\tau_{i+1}=\tau_i\langle n_i,
h(n_i)=\sigma(\tau_i\langle n_i\rangle), h\restriction n_i+1\rangle.$ \qed
\enddemo
\enddemo

Thus for every function $g\in\oo$ the forcing Borel$(\Pi_ng(n))/I_{ie}(g)$ has a dense set consisting of the
$g$-thick trees. It follows easily that the forcing satisfies Axiom A. The following lemma shows that 
$\min\{${\tt cov}$(I_{ie}(g)$ restricted to an arbitrary Borel positive
set), $g\in\oo\}=\min\{${\tt cov}$(I_{ie}(g)):g\in\oo\},$ which will be used in the proof of the relevant variation of (*) in Subsection 2.1.

\proclaim {3.3.5. Lemma}
For every function $g\in\oo$ and every $g$-thick tree $T$ there is a function
$h$ and a continuous map $F:\Pi_nh(n)\to [T]$ such that the preimages of $I_{ie}(g)$-small
sets are
$I_{ie}(h)$-small.
\endproclaim

\demo {Proof}
Fix a function $g\in\oo$ and a $g$-thick tree $T.$ By induction on $n\in\omega$ construct
finite sets $X_n\subset T$ so that $X_0=\{0\}$, for each node $t\in X_n$ there is an
integer $k$ such that the set $X_{n+1}(t)=\{s\in X_{n+1}:t\subset s\}$ consists of
sequences of length $k+1$ and for every $m\in g(k)$ there is a unique $s\in X_{n+1}$ with
$s(k)=m.$ Moreover make sure that $X_{n+1}=\bigcup_{t\in X_n}X_{n+1}(t)$. This is not hard
to do; the sequences in any of the sets $X_n$ will be pairwise incompatible and the union
in the last sentence will always be a union of disjoint sets.

It will be convenient to define the function $h$ so that its range consists of finite sets
rather than natural numbers. Simply let $h(n)=\{Y\subset X_{n+1}:\forall t\in X_n\ |Y\cap
X_{n+1}(t)|=1\}.$ The map $F:\Pi_nh(n)\to [T]$ will be defined by $F(f)=$ the unique
function $e\in [T]$ such that for all numbers $n$ the set $f(n)$ contains an initial
segment of $e$. It is not hard to check the required properties for the function $F.$ \qed
\enddemo

To compare the forcing $PT_g$ of \cite {B1} with the forcing Borel$(\Pi_ng(n))$ modulo $I_{ie}(g)$ note that the former
is a somewhere dense subset of the latter. A moment's thought will then reveal that a suitable iteration of the $PT_g$ forcings
must isolate the invariant {\tt non}(strong measure zero) as well. 

\subhead{3.4. The distributivity of the algebra Power($\omega$) modulo finite}\endsubhead

It is well known that an iteration of Mathias forcing will increase the invariant $\frak h$
defined as the minimum cardinality of a collection of open dense subsets of the algebra $\P(\omega)$ modulo finite
with empty intersection. Actually $\frak h$ is isolated through this iteration. The proof of this fact is a little different 
from the previous cases since Mathias forcing cannot be written as $\P(\R)$ modulo a Borel
generated ideal under any determinacy hypothesis.
It is necessary to settle for a more complicated representation of the forcing. First, some notation.
For sets $a, b\subset\omega$ let $a\subset^*b$ mean that $a$ is included in $b$ up to a finite number of elements.
$[a]$ then denotes the equivalence class of the set $a$ in the algebra $\P(\omega)$ modulo finite, for
a set $A\subset\P(\omega)$ write $[A]=\{[a]:a\in A\}$ and let $I_M$ be the $\sigma$-ideal on $\P(\omega)$
consisting of those sets $A$ for which $[A]$ is nowhere dense in the algebra $\P(\omega)$ modulo finite.

\proclaim {Lemma 3.4.1}
\roster
\item (ZF+DC+AD$\R$) Let $T\subset(2\times Ord)^{<\omega}$ be a tree. Then either $p[T]\in I_M$ 
or there is a condition $p\in\M$ so that $p\Vdash$ the generic real is in $p[\check T]$.
\item (ZFC) Let $T\subset(2\times Ord)^{<\omega}$ be a $<\delta$-weakly 
homogeneous tree, where $\delta$ is a supremum of $\omega$ Woodin cardinals. Then the same dichotomy as in
(1) holds.
\endroster
\endproclaim

\noindent With some additional work, (1) could be restated to say that under ZF+DC+AD$\R$ Mathias forcing is naturally
forcing isomorphic to the algebra $\P\P(\omega)$ modulo the ideal $I_M.$ It is methodologically important to
observe that $I_M$ is not a Borel generated ideal.

\demo {Proof of Lemma 3.4.1(2)}
The following well known geometric condition for Mathias genericity will be used:

\proclaim {Claim 3.4.2}
\cite {SS} Suppose that $a\subset\omega$ is an external $V$-generic Mathias real and $b\subset^* a$ is an infinite external set.
Then $b$ is a $V$-generic Mathias real.
\endproclaim

Now let $T$ be a suitably weakly homogeneous tree such that $p[T]\notin I_M.$ Then there is an infinite set $c\subset\omega$
such that the set $[p[T]]$ is dense below $[c]$ in the algebra $\P(\omega)$ modulo finite. Let $a\subset c$ be a $V$-generic
Mathias real. By the weakly homogeneous absoluteness there is an infinite set $b\in V[a]$ such that $b\in p[T]$ and
$b\subset^* a.$ By the above claim, the set $b$ is a $V$-generic Mathias real and by a wellfoundedness argument involving
the tree $T$ $V[b]\models b\in p[T].$ So there must be a condition $p\in\M$ such that $p\Vdash$ the generic real is in
$p[\check T].$

On the other hand, suppose that some condition $p\in\M$ forces the generic real into $p[\check T]$. Choose a countable
elementary submodel $M$ of a large enough structure containing all the relevant objects and consider the set $A=\{ a\subset\omega:
a$ is a Mathias $M$-generic real meeting the condition $p\}.$ This set is nonempty, Borel by Claim 1.1.3 and its projection
into the algebra $\P(\omega)$ modulo finite is open by Claim 3.4.2. By the choice of the condition $p$ we also have $A\subset p[T].$
Lemma 3.4.1(2) follows. \qed
\enddemo

The proof of the previous lemma also yields

\proclaim {Claim 3.4.3}
Every suitably weakly homogeneous set not in $I_M$ has a Borel subset not in $I_M.$
\endproclaim

With the above facts in hand, the geometric analysis of countable iterations of Mathias forcing proceeds just as in subsection 1.1
replacing the countable ideal by the ideal $I_M$ everywhere, and with the Sacks condition in Definition 1.1.1
replaced by the Mathias condition--splitting into an $I_M$-positive set. The reader is urged to use Lemma 3.4.1(2) to prove on his own 
that the Mathias forcing is forcing isomorphic to the algebra Borel$(\P(\omega))$ modulo $I_M.$ The dichotomy 1.2.1 must be reformulated.
Let $\alpha\in\omega_1$ and let $A\subset\R^\alpha$ be a projective set. Consider the game $G_A$ of $\alpha$
many rounds where at round $\beta\in\alpha$ Eve plays an infinite set $t_\beta\subset\omega$, Adam plays an infinite set $s_\beta
\subset^*t_\beta$ and Eve plays a set $r_\beta\subset^* t_\beta$ in this order. Eve wins if the sequence $\langle
r_\beta:\beta\in\alpha\rangle$ belongs to the set $A.$

Under the assumption of proper class many Woodin cardinals (actually $\omega_1$ many suffice) the game is determined and there are two possibilities.

\roster
\item Adam has a winning strategy
\item Eve has a weakly homogeneous winning strategy and then by an argument essentially identical to that in Subsection 1.2
using Lemma 3.4.3 there is a condition $p\in\M_\alpha$ with $p\subset A.$
\endroster

Note that we could not use a game similar to the original one because there it is important that Adam can play arbitrarily large sets
in the relevant ideal. Here the ideal is not Borel generated and so we would not get a real game and the determinacy of the
game would be open to question.

The argument for the $\frak h$ version of Theorem 0.2 then proceeds exactly as in Section 2 except that the proof
of (*) in Subsection 2.1 has to be changed. Let me recall the setup there. There is a tame invariant $\frak x=\min\{|A|:A\subset\R, 
\phi(A)\land\psi(A)\}$ where the quantifiers of $\phi(A)$ is are restricted to the set $A$ and the natural numbers and 
$\psi(A)=\forall x\in\R\ \exists y\in A\ \theta(x, y)$ where $\theta$ is a formula whose quantifiers range over natural and 
real numbers only. We work in a model where $\frak x<\frak h$ and $A$ is a witness for it, that is $|A|<\frak h,\phi(A)\land\psi(A),$
also we have $\alpha\in\omega_1$ and a condition $p\in\M_\alpha$ and a Borel function $\dot x:p\to\R.$ We want to show
that $\exists y\in A\exists q\leq p\forall\vec r\in q\ \theta(\dot x(\vec r),y).$ 

For each real $y\in A$ let $B_y=\{\vec r\in p:\theta(\dot x(\vec r),y)\}.$ If Eve had a weakly homogeneous winning strategy for one of the games 
$G_{B_y}$ then by (2) above for that real $y\in A$ there would be a condition $q\leq p$ such that $q\subset B_y$ and we would be done.
So it is enough to derive a contradiction from the assumption that Adam has a winning startegy $\sigma_y$ for every game
$G_{B_y},y\in A.$ By a simultaneous transfinite induction on $\beta\in\alpha$ build partial plays of games $G_{B_y}$ for all $y\in Y$
played according to the strategies $\sigma_y$ so that

\roster
\item"{$\circledast$}" writing $t_{\beta y}, s_{\beta y},r_{\beta y}$ for the moves at the $\beta$-th round of the partial play we build for $y\in A$,
the set $r_{\beta y}\subset\omega$ does not depend on $y.$ We can write $r_\beta$ to denote this set.
\item"{$\circledast$}" $\langle r_\gamma:\gamma\in\beta\rangle\in p\restriction\beta$.
\endroster

\noindent To find the moves $t_{\beta y}, s_{\beta y},r_{\beta}$ under the assumption that the partial plays 
$\langle t_{\gamma y}, s_{\gamma y}, r_{\gamma}:\gamma\in\beta\rangle$ were constructed for all $y\in A,$ let $D_y=\{
r\subset\omega:$ there are sets $t,s\subset\omega$ such that the sequence 
$\langle t_{\gamma y}, s_{\gamma y},r_{\gamma}:\gamma\in\beta\rangle\langle t,s,r\rangle$ is a legal partial play of $G_{B_y}$
observing the strategy $\sigma_y\}$. It follows from the definitions that the sets $D_y, y\in A$ are closed under finite changes of
their elements and that they are all open dense in the algebra $\P(\omega)$ modulo finite. Since $|A|<\frak h,$
the intersection of all of these sets is open dense as well and has an element $r_\beta$ in common with the somewhere dense
set $\{r\subset\omega:\langle r_\gamma:\gamma\in\beta\rangle\langle r\rangle\in p\restriction\beta +1\}$. The induction step
is concluded by finding sets $t_{\beta y},s_{\beta y}\subset\omega$ witnessing that $r_\beta\in D_y,$ for all $y\in A.$

Now look at the sequence $\vec r=\langle r_\beta:\beta\in\alpha\rangle.$ Since the strategies $\sigma_y$ were winning for
Adam and in the previous paragraph we produced plays following these strategies whose outcome was the sequence $\vec r,$
it must be that $\theta(\dot x(\vec r), y)$ fails for every real $y\in A.$ This contradicts the property $\psi$ of the set $A.$

\head {4. The tower number}\endhead

The proof of Theorem 0.5 is really just a variation on an argument of Woodin concerning the maximization of
$\Sigma_2$ theory of the model $\langle H_{\aleph_2},\in,\omega_1\rangle.$ Let $\delta$ be a measurable Woodin
cardinal such that for every tame invariant $\frak x,$ if there is a forcing extension satisfying $\aleph_1=\frak x
<\frak t$ then there is such an extension of size $<\delta$. Without loss of generality assume that $2^\delta=\delta^+$.
Let $P_{\frak t}$ be a partial order with the following definition. It is a two step iteration $P_0*\dot P_1$ where $P_0$ is again a
two step iteration $R_0*\dot R_1.$ Here $R_0$ is just the Levy collapse of $\delta$ to $\omega_1$ and $R_1=\{\langle c, D\rangle:
c\subset\delta$ is a closed bounded set such that every limit point $\kappa$ of it is a weakly compact cardinal
of $V$ and the set $c\cap\kappa$ diagonalizes the weakly compact filter on $\kappa.$ The set $D\subset\delta$ belongs
to the weakly compact filter on $\delta$ as computed in $V\}.$ The poset $R_1$ is ordered by 
$\langle c_1, D_1\rangle\leq\langle c_0, D_0\rangle$ if $c_1$ end-extends $c_1,$ $D_1\subset D_0$ and $c_1\setminus c_0\subset D_0.$
Having defined the poset $P_0,$ $P_1$ is just an arbitrary $\sigma$-centered forcing of size $\delta^+=\aleph_2$ in the model
$V^{P_0}$ making MA($\sigma$-centered) and $\frak c=\aleph_2$ true.

The forcing $P_{\frak t}$ deserves an aside. The first step in the iteration defining the poset $P_0$ collapses $\delta$ to $\aleph_1$
and the second step adds a rather mysterious club subset of $\delta$ without adding reals or collapsing $\delta=\aleph_1$ or
$\delta^+=\aleph_2.$ In some sense the model $V^{P_0}$ is supposed to be the most generic model of $\lozenge$ and the model
$V^{P_0*P_1}$ should be the most generic model of MA($\sigma$-centered)+$\frak c=\aleph_2$. The key properties of the forcing
$P_0$ are summed up in the following lemma.

\proclaim {4.1. Lemma} 
\cite {Z} (Woodin) Let $\gamma$ be any Woodin cardinal above $\delta.$
For every poset $Q$ of size less than $\delta$ there are external $V$-generic filters $G_0*G_1\subset
Q*\dot\Bbb P_{<\gamma}$ and $H_0\subset P_0$ so that $G_1\cap\dot\Bbb Q_{<\delta}$ is a $V[G_0]$-generic filter,
$V[G_0][G_2]\subset V[H_0]\subset V[G_0][G_1]$ and $V[G_0][G_1]\models \delta=\aleph_1=|(\P(\delta^+))^V|$.
\endproclaim

Let $j:V[G_0]\to M$ and $i:V[G_0]\to N$ be the elementary embeddings derived from the filters $G_2, G_1$ respectively.
There is a natural factor embedding $k:M\to N$ such that $i=j\circ k$ and since $\delta=\aleph_1^M=\aleph_1^N$, 
necessarily $crit(k)>\delta$. The point in the definition of the forcing $P_0$ is that the model $V[H_0]$ can
be sandwiched between the elementarily equivalent models $M$ and $N$ as far as subsets of $\omega_1$ are concerned.

Back to the proof of Theorem 0.5, let $\frak x=\min\{|A|:A\subset\R, \phi(A)\land\psi(A)\}$ be a tame invariant and let 
$Q$ be a forcing of size less than $\delta$ such that $Q\Vdash\aleph_1=\frak x<\frak t.$ I must prove that
$P_0*\dot P_1\Vdash\aleph_1=\frak x<\frak t$. Choose a Woodin cardinal $\gamma>\delta$ and find the external
objects $G_0,G_1,G_2,H_0$ as in Lemma 4.1 and write $j:V[G_0]\to M$ and $i:V[G_0]\to N$  for the elementary embeddings 
derived from the filters $G_2, G_1$ respectively. The model $N$ is elementarily equivalent to $V[G_0]$ and so
$\frak t>\aleph_1, \frak p>\aleph_1$ and by a theorem of Bell \cite {Be} $MA_{\aleph_1}(\sigma$-centered) are
all true there. Moreover, since the model is closed under $<\gamma$ sequences in $V[G_0][G_1],$ the set
$\P(\delta^{+V})\cap V[H_0]$ is in $N$ and has size $\aleph_1$ there. In particular, the forcing $P_1\in V[H_0]$
is in $N$, it is $\sigma$-centered and by an application of Martin Axiom there is a $V[H_0]$-generic filter
$H_1\subset P_1$ in $N.$ I will show that $V[H_0][H_1]\models\frak x=\aleph_1$ and that will complete
the proof of Theorem 0.5.

Let $A\subset\R\cap V[G_0]$ be a witness to $\frak x=\aleph_1$ in $V[G_0],$ that is $V[G_0]\models |A|=\aleph_1,
\phi(A)\land\psi(A).$ Look at the set $jA\in M.$ Since the set has size $\aleph_1$ in the model and the
critical point of the factor embedding $k:M\to N$ is above $\aleph_1^M=\delta,$ it must be that 
$iA=kjA=jA\in M\subset V[H_0]\subset V[H_0][H_1].$ Now the set $iA$ has the properties $\phi$
and $\psi$ in the model $N$ by elementarity of the embedding $i,$ and it can be argued that it has these
properties in the smaller model $V[H_0][H_1]$ as well. For the property $\phi$ is certainly absolute, and
if $\psi(iA)=\forall x\in\R\exists y\in iA\ \theta(x,y)$ failed in the model $V[H_0][H_1]$ as witnessed
by a real $x$ then it would fail in the model $N$ for the same real $x$, since the reals of both models
are generic extensions of the ground model $V$ and therefore agree on the truth of the projective formula $\theta(x,y).$
Thus $V[H_0][H_1]\models iA$ is a witness to $\aleph_1=\frak x$ as desired.

\enddocument